\documentclass{article}
\usepackage[english]{babel} 
\usepackage{amsfonts}
\usepackage{amsthm}
\usepackage{verbatim}
\usepackage{amsmath}
\usepackage{amssymb}
\usepackage{authblk}
\usepackage[footnotesize]{caption}
\usepackage{graphicx}

\newcommand{\aff}{{\rm aff}}

\newcommand{\config}{{\rm config}}
\newcommand{\tense}{{\rm tense}}

\newcommand{\relcomp}{{\rm relcomp }}
\newcommand{\relbd}{{\rm relbd }}

\newtheorem{thm} {Theorem} [section]
 
\newtheorem{prop}[thm]{Prop.}

\newtheorem{cor} {Corollary} [thm]
\theoremstyle{definition}
\newtheorem{defn}{Definition}
\newtheorem*{rem}{Remark}
\newtheorem{ex}{Example}

\begin{document}
\title{Taut and Tense Networks}
\author{Robert Dawson \\ Saint Mary's University, Halifax, Nova Scotia, B3H 3C3 Canada\\ rdawson@cs.smu.ca}
\date{}

\maketitle
\abstract{Geometrical constructions using flexible cords have been known since the earliest days of recorded mathematics. In this paper we introduce rigorous definitions for two classes of string networks. A \emph{taut} network is one in which all cords are tight in every possible configuration; a \emph{tense} network has configurations in which one or more cords are not tight, but is externally constrained to avoid such configurations. We show that taut networks compute only affine linear functions and subspaces, whereas tense networks (which are closely related to linkages) can trace any algebraic curve.} 
\section{Introduction}
\begin{figure}
    \centering
   \includegraphics[width=0.45\textwidth]{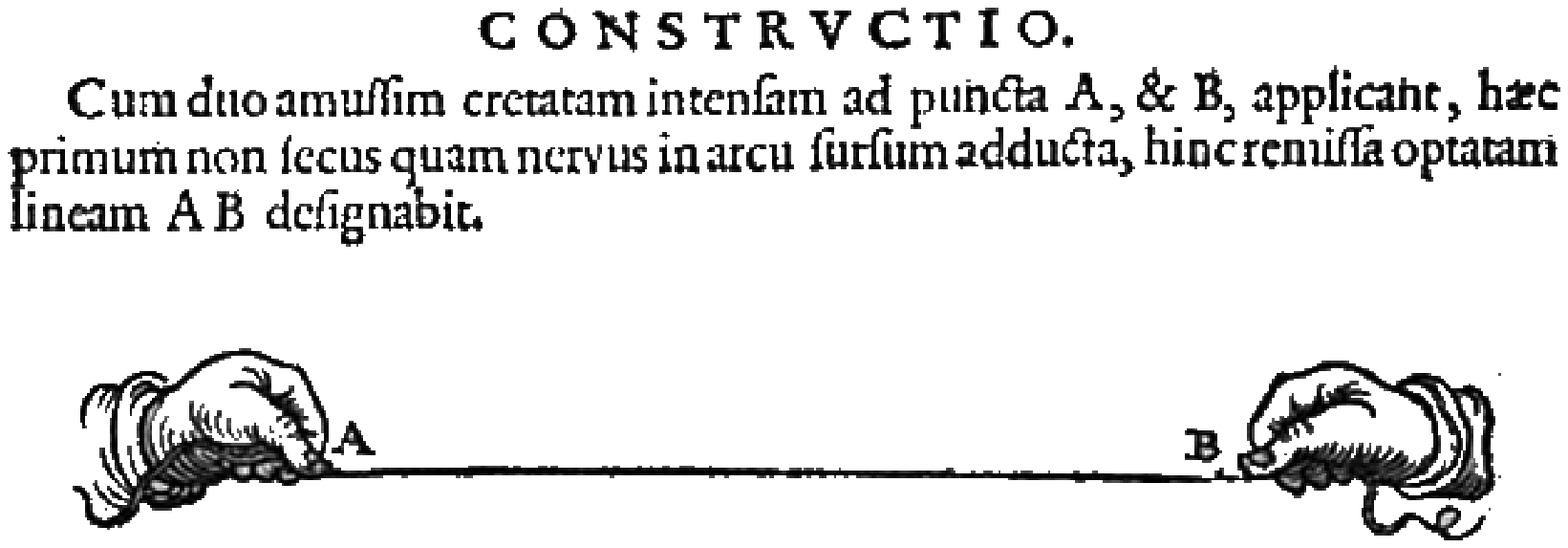}
    \includegraphics[width=0.45\textwidth]{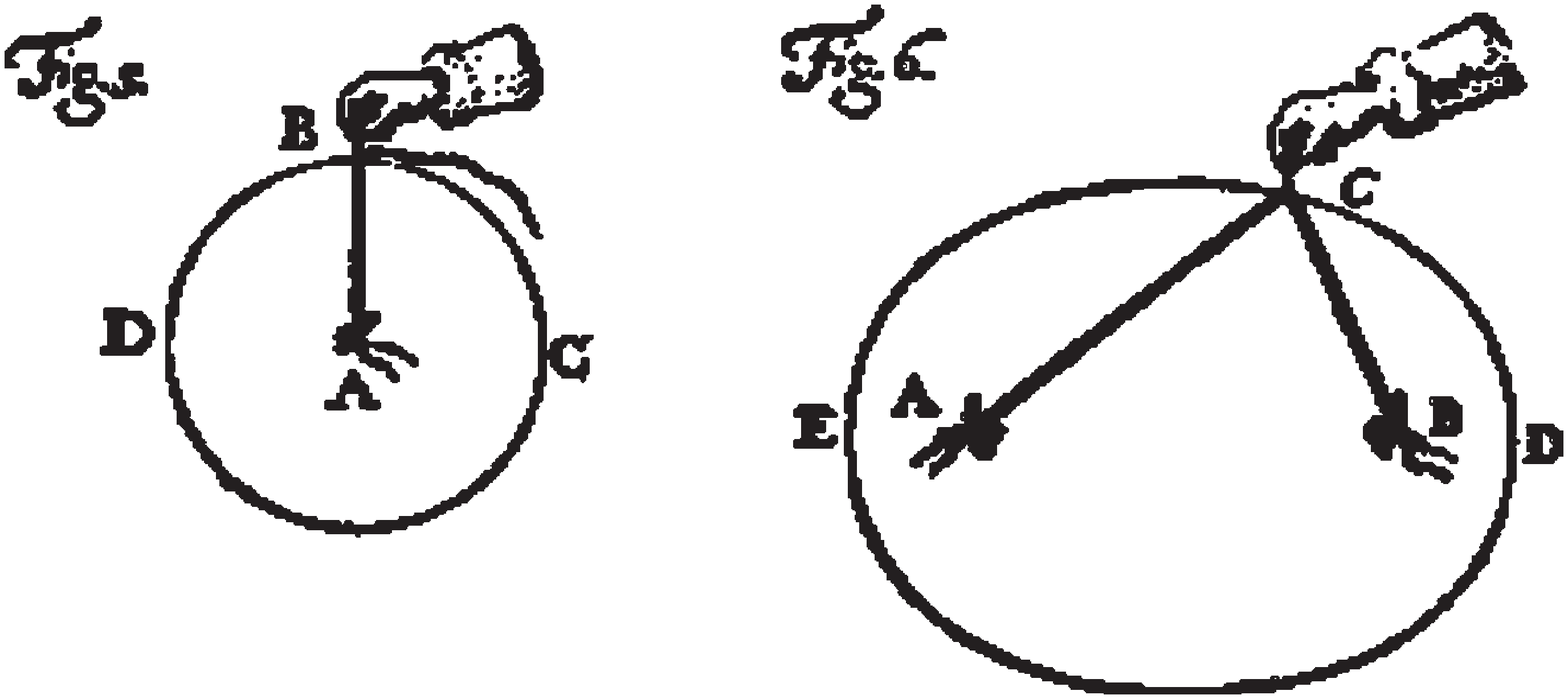}
    \caption{Left: construction of a segment using a taut string, from  Stevin \cite[p. 7]{Stevin}. Right: string-construction of the circle and the ellipse according to von Tschirnhaus \cite[p. 92]{Tschirnhaus}}
    \label{fig:stringhist}
\end{figure}

Constructions using strings under tension have existed for a long time: given the origins of the word ``geometry'' in the Greek words for ``earth measuring,'' it seems likely that a taut cord was the first straightedge and the first compass, as well as the first device for drawing ellipses (Figure \ref{fig:stringhist}). String constructions for Cartesian ovals, hyperbolae, and other shapes have also been given. 

The figure on the left of left of Figure \ref{fig:stringhist}, from the Dutch mathematician Simon Stevin's 1605 \emph{Hymomnematum Mathematicum} \cite{Stevin}, shows a string used as a straightedge. Stevin could have kept it taut with pins rather than hands; the fixed-length cord would in principle stay straight. In contrast, the string compass and ``gardener's ellipse'' from Ehrenfried W. von Tschirnhaus's \emph{Medicina Mentis} \cite{Tschirnhaus} require a continuous outward force to keep the pencil on the proper locus. We will examine the theory behind this apparently minor distinction in Section \ref{Taut}, where we will show that the self-tensioned (or ``taut'') networks exemplified by Stevin's straightedge can only compute affine linear functions, while externally tensioned networks can draw general algebraic curves.

While the compass and ellipse mechanism require external tensioning, they are fairly robust; with a reasonably inelastic cord and moderate force, the curve will be quite accurate.  The stretched line, while simpler, deforms significantly with trivial force. If you try to trace along it with a pencil, even the force needed to keep the pencil against the line may make the traced curve inaccurate. For this reason, when a string is used today as a ruler for landscaping or construction, it is usually dusted with chalk, stretched tight, and then ``snapped'' to leave a mark on the surface. This ``floppiness'' is common (as will be shown) to all self-tensioned networks, and severely limits their utility in some applications. Practical self-tensioned string devices (such as the Etch A Sketch\texttrademark) are usually hybrids with some rigid elements (see \cite{WikiEtch}.)

The string compass can be replaced by a rigid link, transforming it into the simplest example of a mechanical linkage that traces a curve. As we will see below, the classes of networks and of linkages, suitably defined, are equivalent under a bisimulation, and thus can compute the same class of functions.

This paper relates to work in progress with Davide Crippa and Pietro Milici. I would like to acknowledge here helpful discussions with both of them.   

\section{Definitions}

An  \emph{abstract cord} $C$ over a set $N$ consists of a positive real number $L_C$ (the \emph{length} of the cord)  and a list  $(c_0,c_1,...,c_m)$ ($m\geq 1$) of members of $N$ (``nodes'') in which elements may appear more than once but (without loss of any useful generality) not consecutively. (In graph theoretical terms, this is a graph homomorphism $C:P_m\rightarrow K_n$.)  A cord in which the list is a pair will be called a \emph{tie}.

To obtain a general model for a network in $R^d$, we need only add the concept of an anchored node $a$ associated with a location. 

An  \emph{abstract network} $\mathcal{N}$ consists of 
\begin{itemize}
  \item a set $N$ of nodes;
  \item a subset $A\subseteq N$ of \emph{anchored} nodes and a function $\alpha:A\rightarrow R^d$;
  \item a set $\mathcal{C}$ of abstract cords with nodes in $N$.
\end{itemize}

We will call an abstract network \emph{nontrivial} if it has at least one anchored node, at least one non-anchored node, and if every node is on at least one cord.

A \emph{configuration} of an abstract network is a function $\rho:N\rightarrow \mathbb{R}^d$ that restricts on $A$ to $\alpha$, and such that 
\begin{equation}\label{realdef}
d(\rho(c_0),\rho(c_1))+...+d(\rho(c_{n(i)-1}),\rho(c_{n(i)})) \leq L_C
\end{equation} 
for each cord $C\in \mathcal{C}$. A cord is \emph{tense} in a configuration $\rho$ if we have $=$ instead of $\leq$ in \eqref{realdef}, otherwise \emph{slack}; we call $\rho$ itself tense if all its cords are tense.

 For an $n$-node network $\mathcal{N}$ in $\mathbb{R}^d$, the \emph{configuration space} $\config(\mathcal{N})$ is the subspace of $\mathbb{R}^{dn}$ consisting of all configurations; the \emph{tense configuration space} $\tense(\mathcal{N})$ consists of all tense configurations. As just observed, either of these may be empty.  We observe that if a cord is slack in every configuration of a network, removing it does not change the configuration space. We may thus, without loss of generality, omit such cords and assume that every cord of a network is tense in some configuration.  Two nodes that appear as $c_{j-1}$ and $c_j$ on some cord will be called \emph{adjacent.} 

\begin{figure}
    \centering
    \includegraphics[width=0.40\textwidth]{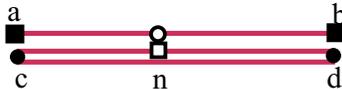}
    \caption{The node $n$ slides on the line $ab$, and is fixed to the loop $ndcn$ that runs around pulleys at $c$ and $d$.}
    \label{fig:nodelegend}
\end{figure}

The constraint \eqref{realdef} models a physical cord with fixed nodes at the ends, and possibly one or more sliding nodes in the middle. A node may be fixed or sliding (on the cord), and, independently, anchored or mobile (in space). Location and being anchored (or not) are properties of the node; if a node appears multiple times on the same or different cords, it necessarily has the same location and anchoring status each time. However, a node may be fixed on one cord and sliding on another. Figure \ref{fig:nodelegend} shows the graphical notation that we will use to represent this: we show a node as square if it is an endpoint (fixed to the cord) and solid if it is an anchor (anchored to the plane.) Note that there are no constraints to prevent nodes, adjacent or otherwise, from occupying the same point in $\mathbb{R}^d$.
 
A physical network may contain (for instance) a string with three fixed nodes, a Y-shaped string,  or a loop with only sliding nodes. Our model loses no generality with respect to its ability to represent the configuration spaces of such networks faithfully. 

 We may consider any point where a string branches as a fixed node, as a sliding node cannot move past it; with this assumption, the complement of the fixed nodes within any network of strings is a disjoint union of real 1-manifolds, which are each either homeomorphic to the interval $I$ or the circle $S^1$. Next, we invoke the ``L\={\i}l\={a}vat\={\i} principle''\footnote{The best-known problem in  Bh\={a}skara II's classic algebra text  \emph{L\={\i}l\={a}vat\={\i}} \cite{Lila} involves a string of pearls broken during lovemaking, the pearls consequently scattering.} that if an interval section does not have a fixed node at each end, it does not constrain the location of any sliding node on it, and may be removed from the network with no effect on the configuration space. An interval section with fixed nodes at both ends is a cord as described above. Finally, a circular loop without fixed nodes imposes no constraint on the configuration space if it has fewer than two sliding nodes. If it does have two or more sliding nodes, one may be replaced by a fixed node without affecting the configuration space; the string is now a cord with the same node at both ends.\\

Our motivating example is the ``gardener's ellipse'' (Figure \ref{fig:stringhist}, right). The tense configuration space can be identified with the ellipse; the configuration space can be identified with the union of the ellipse and its interior. An abstract network may have no configuration: suppose the nodes $a$ and $b$, anchored in the plane at points $(0,-1)$ and $(0,1)$, are the endpoints of a tie of length 1. It may also have a configuration but no tense configuration: for the same nodes, let the tie-length be 3.

If $\mathcal{N}$ is nontrivial, the affine span of $\config(\mathcal{N})$ is never all of $\mathbb{R}^{dn}$ because of its anchored points; thus the concept of metric boundary in $R^{dn}$ is not useful. However, it is convex, so the \emph{affine hull} $\aff(\config(\mathcal{N}))$ will be well-behaved.  Within the affine hull, the \emph{relative interior} of  $\config(\mathcal{N})$ is full-dimensional. Because of the cords, $\config(\mathcal{N})$ is bounded; thus its \emph{relative complement} $\relcomp(\config(\mathcal{N}))$and \emph{relative boundary} $\relbd(\config(\mathcal{N}))$ are also nonempty.

\begin{prop}\label{PropConvex} If $\rho, \sigma \in \config(\mathcal{N})$ , then for $0<t<1$, $\tau := (1-t)\rho+t\sigma \in \config(\mathcal{N})$; and a cord $C$ is tense in $\tau$ and only if it is tense in $\rho$ and in $\sigma$, and $\rho(C_j)-\rho(C_{j-1}) || \sigma(C_j)-\sigma(C_{j-1})$ for each $j$.
\end{prop}
\begin{proof} 
If $a\in A$ then $\rho(a)=\sigma(a)=\alpha(a)$.\\
Suppose that the network has a cord $(x_0,x_1,\ldots,x_n)$ of length $L$. By the triangle inequality, for any two vectors $x$ and $y$, and for $0<t<1$, we have
\begin{equation}\label{TI}
|(1-t)y_0 + ty_1| \leq (1-t)|y_0| + t|y_1|.
\end{equation}
Then
\begin{align*}\label{ConvexCalc}
&\sum_{j=1}^n |((1-t)\rho(x_j)+ t\sigma(x_j))-((1-t)\rho(x_{j-1})+ t\sigma(x_{j-1}))|\\
  = &\sum_{j=1}^n |(1-t)(\rho(x_j)-\rho(x_{j-1}) + t(\sigma(x_j)-\sigma(x_{j-1})|\\
  \leq &\sum_{j=1}^n (1-t)|\rho(x_j)-\rho(x_{j-1})| + \sum_{j=1}^n (1-t)|\sigma(x_j)-\sigma(x_{j-1})|\\ 
  = &(1-t)\sum_{j=1}^n |\rho(x_j)-\rho(x_{j-1})|+ t \sum_{j=1}^n |\sigma(x_j)-\sigma(x_{j-1})|\\
  = &(1-t)L + tL = L;
\end{align*}
and this is an equality if and only if $\rho(x_j)-\rho(x_{j-1}) || \sigma(x_j)-\sigma(x_{j-1})$ for each $j$. 
\end{proof}

\begin{cor}\label{RealConvex}
The configuration space $\config(\mathcal{N})$ is convex in $R^{dn}$, and \emph{a fortiori} path-connected; in particular, a network can be continuously deformed between any two configurations.
\end{cor}  
A network that has only one configuration will be called \emph{static}. For networks, this global property is locally determined. We note that this is not true in general for linkages, tensegrity structures, or for the tense configuration space of a network.  

\begin{cor}\label{Slack} For any given subset $\{C_1,\ldots,C_n\}$ of the cords in a network, if for each cord $C_i$  there exists $\rho_i \in \config(\mathcal{N})$ in which that cord is slack, then for $t_1,\ldots,t_n > 0$, $t_1+\cdots+t_n=1$, we have that the $C_i$ are simultaneously slack in the configuration $t_1\rho_1+\cdots+t_n\rho_n \in \config(\mathcal{N})$. 
\end{cor}

Linear interpolations between two configurations are always themselves configurations. However, if a linear \emph{extrapolation} of two tense configurations is itself a configuration, this has important implications.
\begin{prop}\label{Extrapolation} 
  If $\rho, \sigma \in \config(\mathcal{N})$, if cord $C$ is tense in $\rho$ and $\sigma$, if $t<0$ or $t>1$, and if $\tau := t\rho + (1-t)\sigma \in \config(\mathcal{N})$, then
\begin{enumerate}
  \item $C$ is tense in $\tau.$ 
  \item If $\upsilon := s\rho +(1-s)\sigma \in \config(\mathcal{N})$, then  $C$ is tense in $\upsilon$. 
\end{enumerate}
\end{prop}
\begin{proof}(1) Suppose (for a contradiction) that $C$ is slack in $\tau$ and (WLOG) that $t<0$. Then by the previous proposition $C$ would be slack in $\rho$, contradicting our hypothesis.\\  
(2) It only remains to prove the result for $0<s<1$; this follows from:
\begin{align*}
\rho &= \frac{\tau}{t} -\frac{1-t}{t}\sigma\\
s\rho +(1-s)\sigma &= \frac{s}{t}\tau -\frac{s(1-t)}{t}\sigma +(1-s)\sigma\\
&= \frac{s}{t}\tau + \left (1 - \frac{s}{t}\right) \sigma.\\
\end{align*}
As $t<0<s$, we have $s/t < 0$, and by part 1, $s\rho +(1-s)\sigma$ is tense.
\end{proof}

The ``gardener's ellipse'' shows that we need the assumption that the third configuration is not between the two given tense configurations.

\section{Taut Networks}\label{Taut}

In this section we will define and study taut networks: these are the ``pretensioned'' networks mentioned above, for which every configuration is tense. We will show that a curve or function can be computed by a taut network if and only if it is affine linear; and we will show that all taut networks are subject to the type of low-force deformation observed in the chalkline. We begin with a proof that a cord which is locally always tense is globally always tense.

\begin{prop}\label{prop:LocalGlobal}
If $C$ is a cord of $\mathcal{N}$, $\rho \in \config(\mathcal{N})$, and there exists $\epsilon>0$ such that whenever $\sigma \in \config(\mathcal{N})$ and $|\rho-\sigma|<\epsilon$ then $C$ is tense in $\sigma$, then $C$ is tense in all of $\config(\mathcal{N})$.
\end{prop}
\begin{proof} Let $\tau \in \config(\mathcal{N})$. By Proposition \ref{PropConvex}, $t\tau +(1-t)\rho \in \config(\mathcal{N})$ for $0<t<1$. But, by hypothesis, for small enough $t>0$, $C$ is tense in that configuration; and by Proposition \ref{Extrapolation} this extends to all $t\in[0,1]$.
\end{proof}

 A network is defined to be \emph{taut} if every cord has this property. The following is an immediate consequence of Proposition \ref{PropConvex}:

\begin{prop}\label{PropPar} 
If $\rho,\sigma \in \config(\mathcal{N})$, and $C$ is a cord of $\mathcal{N}$, then $\rho(C_j)-\rho(C_{j-1}) || \sigma(C_j)-\sigma(C_{j-1})$ for each adjacent pair of nodes $C_{j-1},C_j$ of $C$.
\end{prop}

We will define a node to be \emph{immovable} if its position is the same in every configuration, regardless of whether it is directly anchored, and \emph{movable} otherwise. 

\begin{figure}
    \centering
    \includegraphics[width=0.75\textwidth]{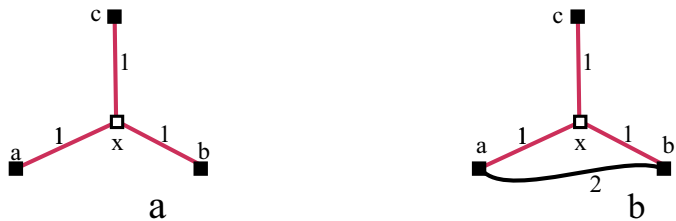}
    \caption{The node $x$ is immobile but not anchored.}
    \label{fig:Ynet}
\end{figure}

\begin{ex} Consider the \emph{Y network} (Figure \ref{fig:Ynet}a) that has $A=\{a,b,c\}$, the set of vertices of an equilateral triangle with edge 1; $x$ is the only other vertex, and $(a,x),(b,x)$, and $(c,x)$ are all cords with length $\frac{\sqrt{3}}{3}$.  Then  Y is tense, static (and vacuously taut), and $x$ is immovable but not anchored.\end{ex}

\begin{figure}
    \centering
    \includegraphics[width=0.6\textwidth]{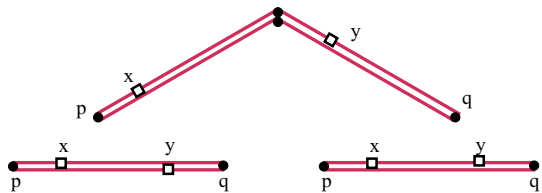}
    \caption{Variations on the clothesline network: the position of $y$ copies that of $x$.}
    \label{fig:Clothesline}
\end{figure}

\begin{ex}
The \emph{clothesline} (figure \ref{fig:Clothesline}) has $N = \{p,x,q,y\}$, $A=\{p,q\}$, and two cords $(x,p,y)$ and $(y,q,x)$ with length $L = d(\alpha(p),\alpha(q))$. Every configuration of this has $x$ and $y$ on the line segment between $\alpha(f)$ and $\alpha(g)$. All such configurations are tense, thus they are all taut, and $d(\rho(x),\alpha(p))= d(\rho(y),\alpha(q))$. This network, and variations on it, may therefore be used as a ``data bus'' to transfer information within networks.
\end{ex}

\begin{ex} Let the Y network be augmented with a cord of length 2 connecting $a$ and $b$ (Figure \ref{fig:Ynet}b.) This network is static, but trivially not tense (and hence not taut.) But clearly this, and any static network, may be made taut by removing nonfunctional cords.
\end{ex}

\begin{figure}
    \centering
    \includegraphics[width=0.40\textwidth]{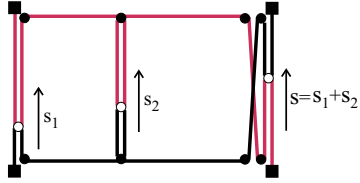}
    \caption{The position of the node $s$ is the sum of the positions of the summand nodes $p_i$.}
    \label{fig:Adder}
\end{figure}

\begin{ex} The \emph{adder} is shown in Figure \ref{fig:Adder}. The vertical position of the node $s$ is the sum of the vertical positions of $s_1$ and $s_2$.
\end{ex}

Other examples of taut networks include the \emph{Cartesian} network and the related \emph{scaler.}
\begin{figure}
    \centering
    \includegraphics[width=0.40\textwidth]{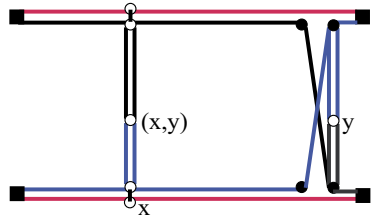}
    \caption{The position of the node $(x,y)$ is the Cartesian product of the positions of $x$ and $y$.}
    \label{fig:Cartesian}
\end{figure}   

\begin{ex}\label{CartEx}
The \emph{Cartesian} network is essentially an Etch A Sketch\texttrademark without rigid components. Versions exist in any dimension $d\geq2$; the plane version is shown in Fig. \ref{fig:Cartesian}. It can do any of the following:
\begin{itemize}
  \item The $d$-dimensional coordinates of the movable node $x$ are given by the movable nodes $x_1,\ldots,x_d$; the network simultaneously computes all the projections.
  \item  The cartesian product of the locations of the movable nodes $x_1,\ldots,x_d$ is given by the movable node $x$.
  \item If other cords constrain $x$ to a subset of $R^d$, this enforces a relation on the movable nodes $x_1,\ldots,x_d$. In particular, if $x$ is constrained to the graph of a function $x_1 = f(x_2,\ldots,x_d)$, the network computes that function. 
\end{itemize} 

\end{ex}

\begin{figure}
    \centering
    \includegraphics[width=\textwidth]{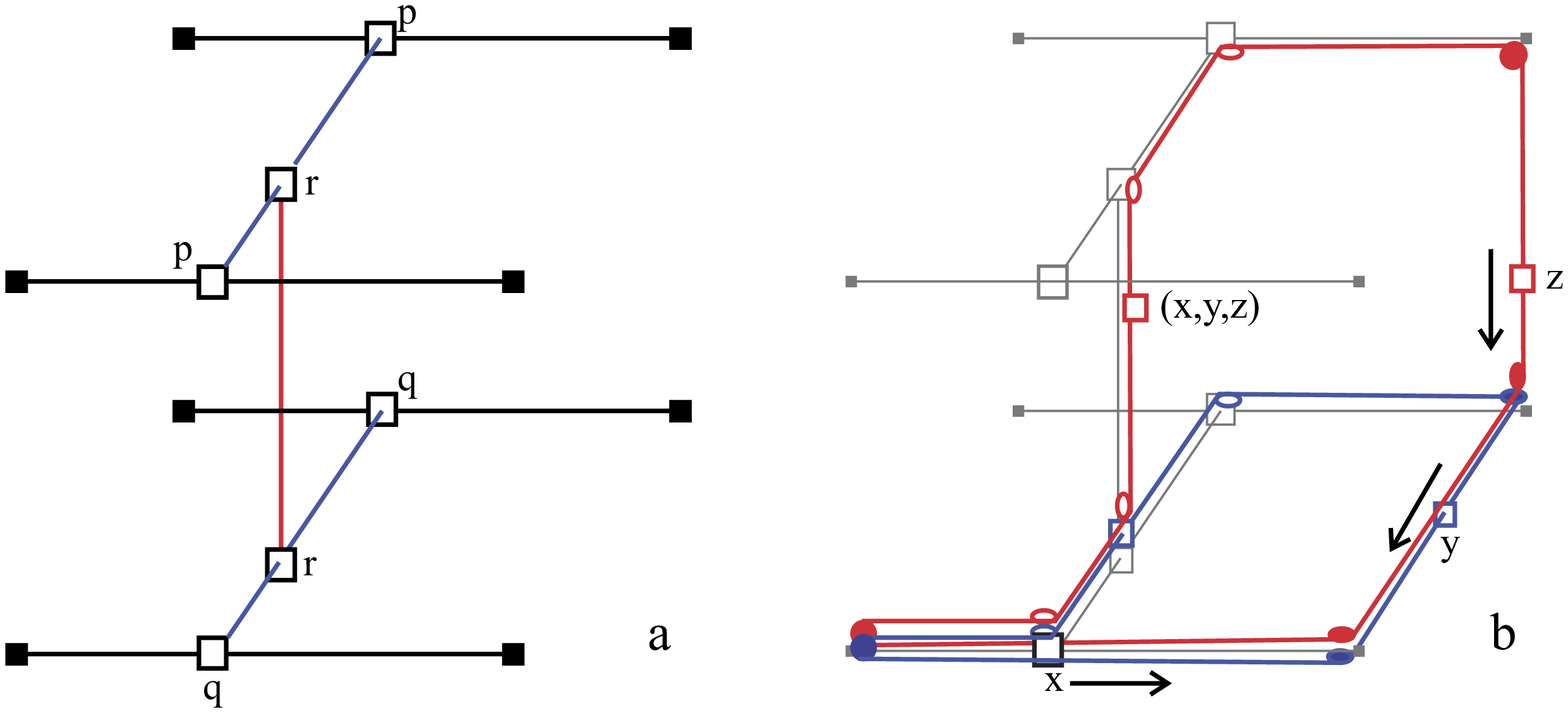}
    \caption{The position of the node $(x,y,z)$ is the Cartesian product of the positions of $x$,$y$, and $z$.}
    \label{fig:Cart3}
\end{figure}

\begin{ex} Fig. \ref{fig:Cart3} shows a three-dimensional version of the Cartesian network. On the left we see the skeleton: the cords $pp'$ and $qq'$ must remain perpendicular to the cords that their endpoints move on. Moreover, $rr'$ not only remains perpendicular to $pp'$ and $qq'$, but keeps them in a plane perpendicular to the four fixed cords, so that all cords shown remain parallel to the coordinate axes. The diagram on the right adds two more loops that copy the $y$ and $z$ coordinates onto intervals between anchored nodes. 
\end{ex}

\begin{figure}
    \centering
    \includegraphics[width=0.40\textwidth]{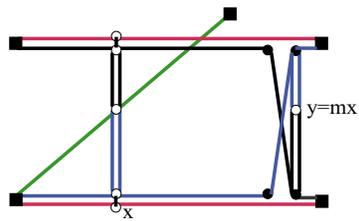}
    \caption{The position of the node $y$ is that of the node $x$ multiplied by a hardwired constant $m$.}
    \label{fig:Scaler}
\end{figure}

\begin{ex} As a simple variation of Example \ref{CartEx}, if $x$ is constrained to lie on a straight line, we can compute $y=mx$ for any constant $m$ (Figure \ref{fig:Scaler}); this network is the \emph{scaler}.
\end{ex}

The constant $k$ is built into the network. It might seem that the network could be modified to make $k$ variable and create a multiplier/divider, but as we will see below this cannot be done in a taut network.

\begin{thm}
The coordinates $\rho_i(n_k)$ of the mobile points of a taut network $\mathcal{N}$ satisfy a linear relation for all $\rho \in \config(\mathcal{N})$.
\end{thm}
\begin{proof}
Suppose the network to consist of $K$ cords $(C_i:1\leq i\leq K)$ , each with nodes $(n(i,j), 0\leq j \leq k(i)$. (As always, nodes may be shared between cords, and repeated on one cord, provided that they are nonadjacent.) As the network is taut, for any two adjacent nodes $n(i,j-1)$ and $n(i,j)$ and any two configurations $\rho,\sigma \in \config(\mathcal{N})$, we have (by Proposition \ref{PropPar}) that
$$\rho(n(i,j))-\rho(n(i,j-1))\;\|\;\sigma(n(i,j))-\sigma(n(i,j-1)).$$
Thus we obtain
\begin{equation}\label{ParEq_Sd}
(\rho(n(i,j))-\rho(n(i,j-1))) = l(\rho,i,j)(\vec{c}_{ij})
\end{equation}
where $\vec{c}_{ij}$ is a unit vector independent of $\rho$. This gives $d$ linear equations for each adjacent pair of nodes.
Summing along the cord we get
\begin{equation}\label{Total_Kd}
\sum_{j=1}^{k(i)} l(\rho,i,j) = L_i.
\end{equation}
In a taut network, an end node of a cord must be connected to something, and other nodes may be.  Such a node may be anchored:
\begin{equation}\label{Anch_d}
\rho(n(i,j)) = \alpha(n(i,j)) ;
\end{equation}
or it may be on another (or the same) cord
\begin{equation}\label{Cord_d}
\rho(n(i,j)) = \rho(n(i',j').
\end{equation}
All the constraints on the motion of a taut network are of one of the four types (\ref{ParEq_Sd}-\ref{Cord_d}), and they are all given by linear equations in the variables $\rho(i,j)$ and $l(\rho,i,j)$. By assumption this system has at least one solution. If it is underdetermined, the network has one or more degrees of freedom, described by a linear relation.
\end{proof}

\begin{cor}There is no taut network that constrains a point to move on a circle or ellipse; and there is no taut network that computes the product of two lengths.
\end{cor}
This observation may be thought of as justifying the hands shown in Tschirnhaus's diagrams (Fig. \ref{fig:stringhist}, center and right.) 

\begin{thm}
Any affine relation on the coordinates of a set P of points in $R^d$ can be realized by a taut network, the nodes of which are a superset of P, using $O(2^d\cdot|P|) segments$.
\end{thm}
\begin{proof}This follows from the existence of the Cartesian, Scaler, Adder, and Clothesline networks. A Cartesian network for each point of P extracts its coordinates $p_{ij}$; each of these is either held constant, left free, or passed to a Scaler network that computes $a_{ij}p_{ij}$; and an Adder constrains $\sum a_{ij}p_{ij}$ to be constant.
\end{proof}

We conclude that taut networks can compute (over a bounded range) \emph{exactly} the affine relations. 

\begin{defn} The possible locations of a specified node $x$ in a taut configuration form a convex subset of $R^d$ (possibly a single point.) We will call the dimension of this set the \emph{mobility} of $x$.
\end{defn} 

\begin{prop}\label{PropMob}
\begin{enumerate}
\item If nodes $x$ and $y$ are adjacent in a taut network, then $m(x)-m(y) \in \{-1,0,1\}$.
\item If a node $x$ in a taut network is adjacent to nodes $y_1,\ldots,y_n$ with lower mobility, then, in any configuration $\rho$, the points $\{\rho(x)$, $\rho(y_1), \ldots,$ $\rho(y_n)\}$ are collinear.
\end{enumerate}
\end{prop}
\begin{proof}
Suppose that $m(x)\neq m(y)$. Then, without loss of generality, suppose that $m(x) > m(y)$. There must be configuration $\rho,\sigma$ such that $\rho(x)\neq \sigma(x)$ but $\rho(y)=\sigma(y)$.
By Proposition \ref{PropConvex}, $\rho(x)-\rho(y)\; \|\; \sigma(x)-\sigma(y)$. Thus we have 
$$\rho(x)-\sigma(x)\; \|\; (\rho(x)-\rho(y)) - (\sigma(x)-\sigma(y))\; \|\; \rho(x)-\rho(y).$$
As this is true for all such $\sigma$, (1) follows. But this is true for any $y$, so 
$$\rho(x)-\rho(y_i)\; \|\; \rho(x)-\sigma(x)\; \|\; \rho(x)-\rho(y_j).$$ 
which gives us (2). 
\end{proof}

\begin{figure}
    \centering
    \includegraphics[width=0.40\textwidth]{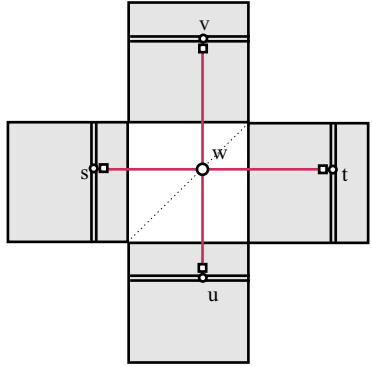}
    \caption{The node $w$ is adjacent only to nodes of higher mobility.}
    \label{fig:CEx}
\end{figure}

\begin{ex}
Consideration of simple cases suggests the conjecture that every movable node in a taut network must be adjacent to a node of equal or lower mobility, but this is in fact not the case. In Figure \ref{fig:CEx}, the points $u,u',v,v'$ are constrained by linked Cartesian networks (shown shaded, with details omitted) to have the following coordinates:
\begin{align*}
\rho(u) &= (x-1,y)\\
\rho(u') &= (x+1,y)\\
\rho(v) &= (y-1,x)\\
\rho(v') &= (y+1,x)
\end{align*}
Cords of constant length 2 join $u$ to $u'$ and $v$ to $v'$, and the point $w$ is at the intersection of those cords. It is clear that $\rho(w) = (y,y)$, so that $m(w)=1$ while 
$m(u)=m(u')=m(v)=m(v')=2$.
\end{ex}
 
\section{Tautness and Firmness}
Our assumption that the length of each cord is preserved is equivalent to an assumption of perfectly inelastic cords. This is good pure mathematics, but poor engineering. Any applied force may result in an increase in the length of one or more cords. In most cases the extension is proportional to the force; but to move a node perpendicularly to a straight cord requires a force of a smaller order of magnitude. In the words of Whewell's 1819 \emph{Elementary Treatise on Mechanics}:

\begin{quote}Hence no force however great can stretch a cord however fine into an horizontal line which is accurately straight: there will always be a bending downwards\footnote{The first 19 words of this have enjoyed some notoriety as a ``found poem'' - see for example Vol. IX of Bliss Carman's \emph{The World's Best Poetry} \cite{Car}} \cite{Whew}.
\end{quote}

It is natural to ask: can a nonvacuously taut network be constructed in such a way that displacing any node by $\epsilon$ from its theoretical locus will require an extension of the same order in at least one cord? (Equivalently, if the cords of a network are replaced by elastic cords each subject to Hooke's law, does a similarly linear force law apply to any departure from the original configuration space?)   

\begin{defn}Let $\mathcal{N}_\epsilon$ be the abstract network derived from $\mathcal{N}$ by increasing each cord length by $\epsilon$; any configuration of $\mathcal{N}_\epsilon$ is also a configuration of $\mathcal{N}$. Given two configurations $\rho, \rho'$ of an abstract network, we define the distance $d(\rho,\rho'):=\max_j(\|\rho(c_j)-\rho'(c_j)\|)$. We will say that $\mathcal{N}$ is \emph{firm} if there exists $k$ such that for all $\epsilon>0$ and all configurations $\rho$ of $\mathcal{N}$ and $\rho'$ of $\mathcal{N}_\epsilon$, the distance $d(\rho,\rho')$ is always less than $k\epsilon$.
\end{defn}

\begin{figure}
    \centering
    \includegraphics[width=0.70\textwidth]{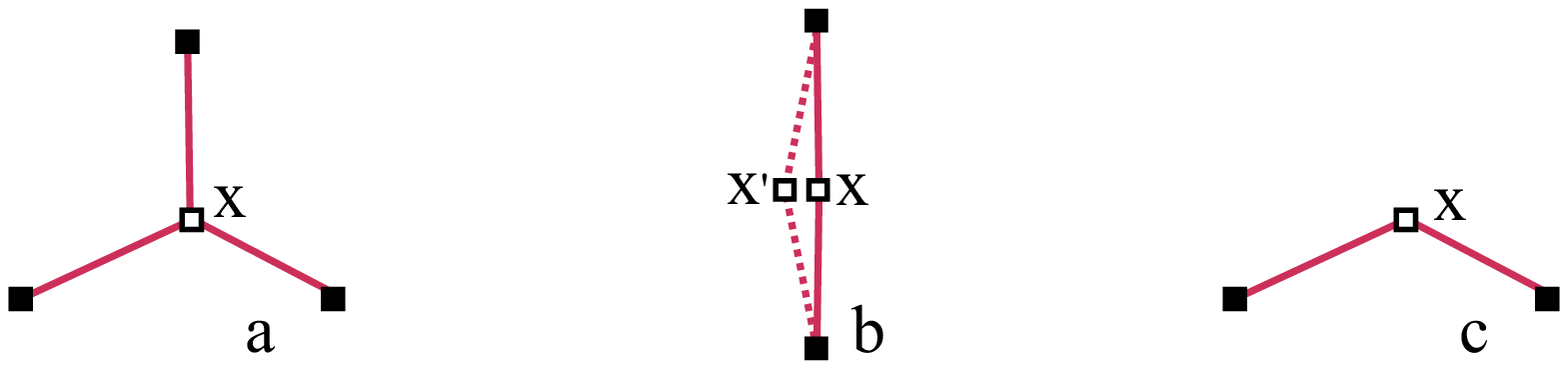}
    \caption{Networks $a$ and $c$ are firm; network $b$ is not.}
    \label{Firm}
\end{figure}

\begin{ex} The vacuously taut network in Figure \ref{Firm}a is firm. The taut network in Figure \ref{Firm}b is not. The non-taut network in Figure \ref{Firm}c is also firm.
\end{ex} 

\begin{prop}
For $d>1$, no nonvacuously taut network in $R^d$ is firm.
\end{prop}

\begin{proof} Consider the set $M$ of movable nodes; if the network is nonvacuously taut, this is nonempty. If no element of $M$ is adjacent to an immovable node, then the network is disconnected, with no immovable node in the same component as any movable node, and we can translate $M$ freely. Suppose, then, that some node $m\in M$ is adjacent to one or more immovable nodes $\{a_1,\ldots,a_n\}$; then by Proposition \ref{PropMob} the nodes $\{m,a_1,\ldots,a_n\}$ all lie on one line. The point $m$ must be free to move (for some distance) along this line; for if some cord stops its motion in one direction, then that cord would lose tension when $m$ was moved in the opposite direction.

Suppose that no cord in the network has more than $N$ sections between nodes, and that the minimum length of any section is $a$. If the lengths of the cords are increased by 
$\epsilon$, each section may be increased by $\epsilon/N$. Then any movable node $m$ may be moved in any direction perpendicular to its associated line by a distance $\sqrt{(a+\epsilon/N)^2-a^2} > \sqrt{(2a\epsilon)/N}$. It follows that it may be moved by at least that amount in any direction whatsoever.

But then we may translate the entire set $M$ a distance $\sqrt{(2a\epsilon)/N}$ in any direction whatsoever, as only the distances between mobile and immobile points changes. This is the distance between the two configurations, and as $\epsilon\rightarrow 0$ we have $d(\rho,\rho')/\epsilon = \sqrt{(2a\epsilon)/N}/\epsilon \rightarrow \infty$, and the network is not firm.
\end{proof}

\section{Conditionally Tense Networks}

A non-taut network may have tense configurations;  we will call such a network ``conditionally tense." The simplest examples are the ``string compass'' and the ``gardener's ellipse'' of Figure \ref{fig:stringhist}  in their working configurations. As shown above, such a tense configuration is always on the relative boundary of the network's configuration space. 

\begin{prop}\label{Tenseprop}Let $\mathcal{N}$ be a non-taut network. If $\rho \in \relbd(\config(\mathcal{N}))$, then $\mathcal{N}$ has at least one cord for which \eqref{realdef} is an equality.
\end{prop}
\begin{proof} If  \eqref{realdef} is not an equality for any cord, then by continuity every nonanchored node can move freely in an $\epsilon$-ball within the configuration space while keeping \eqref{realdef} strict. The product of these balls contains an $\epsilon$-ball of the affine hull of the configuration space, contradicting the assumption that the configuration is on the relative boundary.
\end{proof}

\begin{prop}
Let $\mathcal{N}$ be a non-taut network; then $\tense(\mathcal{N}) \subseteq \relbd(\config(\mathcal{N}))$.
\end{prop}\begin{proof}
Let $\rho \in \tense(\mathcal{N})$. By Proposition \ref{prop:LocalGlobal}, if $\mathcal{N}$ is non-taut, then for any $\epsilon>0$ there exists a non-tense $\sigma\in\config(\mathcal{N})$ with $\|\rho-\sigma\|<\epsilon$. Define $\tau := 2\rho-\sigma$; clearly $\tau \in \aff(\mathcal{N}))$, $\|\rho-\tau\| < \epsilon$, and (by Proposition \ref{Extrapolation}) $\tau\not\in\config(\mathcal{N})$. Thus $\rho \in \relbd(\config(\mathcal{N}))$. 
\end{proof}

\begin{figure}
    \centering
    \includegraphics[width=0.8\textwidth]{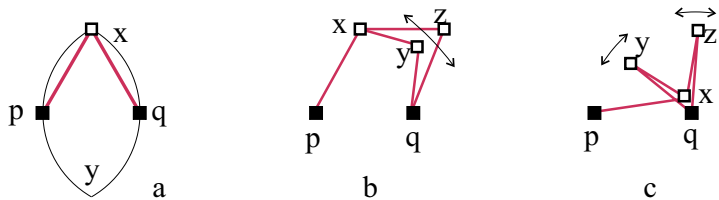}
    \caption{(a) $\tense(\mathcal{N})=\{x,y\}$; $\relbd(\config(\mathcal{N}))$ is the union of the two arcs. (b,c) The number of degrees of freedom varies within the tense configuration space. The nodes $x$ and $q$ coincide in (c).}
    \label{fig:VesPis}
\end{figure}

\begin{ex} The containment may be proper: consider the network of Figure \ref{fig:VesPis}a, for which 
$\tense(\mathcal{N})= \{x,y\}$, while $\relbd(\config(\mathcal{N}))$ is the union of the two circular arcs (the ``vesica piscis'' or ``mandorla''). Note also that neither $\tense(\mathcal{N})$ nor $\relbd(\config(\mathcal{N}))$ is convex in this example, although their convex hulls are (by Proposition \ref{RealConvex}) within $\config(\mathcal{N})$.
\end{ex}

\begin{ex} The local dimension of the tense configuration space of a network may vary. Figure \ref{fig:VesPis}$b$ and $c$ show different configurations of the same network. In (b), the position of node $x$ determines (locally, assuming the cords remain tense) the positions of $y$ and $z$. In (c), where $\rho(x)=\rho(q)$, the nodes $y$ and $z$ are free to move independently. The 1-manifolds of tense configurations resembling (b) connect to the 2-manifold of tense configurations like (c) at the four configurations where all nodes lie on one line.

Note that this cannot happen if the network is taut; convexity requires that the neighborhood of every point must have the same number of degrees of freedom.

\end{ex}

The reader may note that a tense tie is essentially a rigid two-node \emph{link}. A \emph{simple plane linkage} is defined to be a set of nodes $N$ with a specified set of two-node links $\{L_i\}$, each consisting of a pair of nodes $n_i,n'_i \in N$ and a distance $d_i>0$. A realization of a plane linkage is a function $\varrho:N\rightarrow R^2$ such that $\|\varrho(n_i)-\varrho(n'_i)\| = d_i$ for each link, and $\varrho(n_i) = \alpha_i$ if $n_i$ is anchored. 

\begin{figure}
    \centering
    \includegraphics[width=0.8\textwidth]{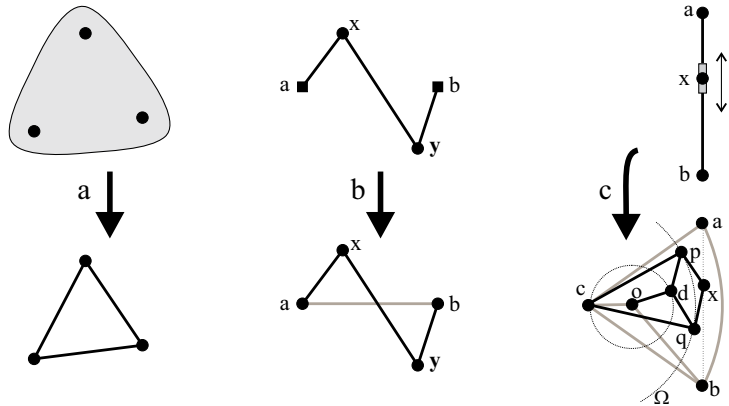}
    \caption{Simulations using simple links of (a) a multi-node link; (b) anchored nodes; (c) a slider.  Grey links in the lower row are ``backplane'' serving only to keep nodes in place.}
    \label{fig:LinkSim}
\end{figure}

Assemblages of simple links can simulate more complicated elements. Among these:
\begin{description}
  \item[Multi-node links:] We could allow a link  $\Lambda_i$ to connect three or more nodes in a specified planar position relative to each other; but any such structure can be replaced (see, for instance, \cite{KM}) by a complete graph of two-node links, or indeed a suitable subgraph: see Figure \ref{fig:LinkSim}a.
  \item[Anchored nodes:] We could specify a set $A$ of anchored nodes $n_i$ each with a specified location $\vec{a}_i$ in the plane; but this can be replaced (up to isometry) by a ``backplane'' in the form of an additional multinode link. (For this reason, the literature often refers to Watt's linkage, a chain of three links with its end nodes connected to fixed shafts, as a ``four-bar linkage'': see Figure \ref{fig:LinkSim}b)
  \item[Sliders:] We can expand our ``parts list'' to include a node sliding freely on a fixed-length rod. This can be simulated using Peaucellier's inverting cell, shown in Figure \ref{fig:LinkSim}c. The quadrilateral $pdqx$ is always a rhombus, and $p,q$ always lie on the same circle $\Omega$ about $c$. Thus $c,d,x$ are collinear, and by similar triangles we have $|cd||cx| = |cp|^2$. The points $d$ and $x$ are thus inverses in $\Omega$, and if the link $od$ constrains $d$ to move on a circle passing through $c$, $x$ must move on a straight line.
\end{description} 
 
\begin{figure}
    \centering
    \includegraphics[width=\textwidth]{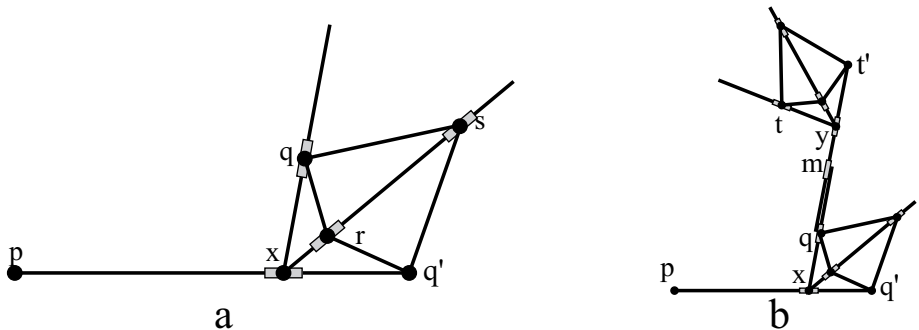}
    \caption{Simulations using simple links and sliders of (a) a cord with one mobile point $x$; (b) a cord with two mobile points $x,y$.}
    \label{fig:Flex}
\end{figure}

\begin{thm} Any algebraic plane curve can be locally traced by a node of some tense plane network; and any curve that can be locally traced by a node of a tense plane network is algebraic.
\end{thm} 
\begin{proof}  Kempe showed \cite{K75} that in a linkage with one degree of freedom, each moving node traces an algebraic curve, and conversely that every algebraic curve may be traced locally by such a linkage. (Kempe's original proof had gaps; a complete proof was given by Kapovich and Millson \cite{KM}, and a more elementary complete proof by Power \cite{Power}.  A good introduction is given in \cite{DO}.) We shall show that linkages and tense networks can model each other.  

Each link of a linkage can be replaced by a tie with the same nodes and length; clearly if that tie is tense, the length is maintained. Thus any linkage may be modeled by a tense network.

To model a tense network by a linkage, we replace every tie by a link. A tense cord with a single movable node is replaced by a linkage as shown in Figure \ref{fig:Flex}a. In this linkage, the length $px$ and the angle $\angle pxq$ are both variable, but $\triangle sqr \cong \triangle sq'r$, $\angle rsq \cong \angle rsq'$, and so $xq = xq'$, and $px+xq = pq'$ which is constant. Thus, for instance, if we anchor $p$ and $q$, the locus of $x$ will be (part of) an ellipse.

A cord with two (or more) links is modeled by a concatenation of such linkages (Figure \ref{fig:Flex}b). Two sliders keep the points $x,q,m,$ and $t'$ collinear, and a link keeps $qt'$ constant. It follows that $px+xy + yt = pq'+qt'$ which is constant, though the positions of $x$ and $y$, and the angles there, can vary.
\end{proof}

\begin{rem} We can consider a third class of networks, in which not all cords are tense in any one configuration. The simplest example is the network of Figure \ref{fig:VesPis}a, which differs from the ``gardener's ellipse'' in that the node $x$ is fixed. As this network moves along the relative boundary of its configuration space, the node $x$ traces the \emph{vesica piscis} curve, which is not given by any one algebraic equation. However, it is piecewise algebraic; we can consider it as the boundary of the intersection of two discs, each the configuration space of a subnetwork in which one cord is omitted. Each of these configuration spaces is convex; therefore their intersection is convex too. 

These observations extend to more complicated networks in which the relative boundary of the configuration space, rather than the tense configuration space, constrains a node to trace a curve. The curve traced by any such network is obviously piecewise algebraic, and this class of networks does not appear to have great theoretical interest.
\end{rem}

\section{Conclusions}
We have given a simple model for geometric calculation devices in Euclidean space using only strings with fixed or sliding nodes. We have seen that such networks divide naturally into two classes: those that maintain their own tension and those requiring an external force to keep them tense. The first class can compute linear functions of any number of variables, and the only one-parameter loci that it can trace are straight line segments; the second have the same computational power as linkages, and can compute algebraic functions and trace (portions of) algebraic curves.

\end{document}